%% file: main.tex
\documentclass{article}

\usepackage{aem}
\usepackage{amssymb}
\usepackage{amsmath}
\usepackage{epsfig}

\title{An introduction to fractional calculus \\ {\small Numerical methods and application to HF dielectric response }}

\name{Andr\'e Persechino}

\address{
LitCOMP - SI: Signal \& Image Laboratory, \\
Brazilian Center Physics Research, Brazil \\
*corresponding author, E-mail: {\tt aamerico@cbpf.br}
}

\begin{document}
\maketitle

\begin{abstract}
    \input{abstract/sec_abstract}
\end{abstract}

\section{Introduction}                          \label{sec_introduction}
    \input{introduction/sec_introduction}

\section{Fractional calculus}                   \label{sec_fractional_calculus}
    \subsection{Fundamental concepts}           \label{sec_fundamental_concepts}
        \input{fractional_calculus/sec_fundamental_concepts}
    \subsection{Riemann-Liouville approach}     \label{sec_RL_approach}
        \input{fractional_calculus/sec_RL_approach}
    \subsection{Gr\"unwald-Letnikov approach}   \label{sec_GL_approach}
        \input{fractional_calculus/sec_GL_approach}
    \subsection{On the plurality of definitions: identification criteria}   \label{sec_identification_criteria}
        \input{fractional_calculus/sec_identification_criteria}
        
\section{Universal dielectric response} \label{sec_universal_dielectric_response}
    \subsection{Linear media and temporal non-locality} \label{sec_linear_media}
        \input{universal_dielectric_response/sec_linear_media}
    \subsection{Empirical models for polarization}  \label{sec_empirical_models}
        \input{universal_dielectric_response/sec_empirical_models}
    \subsection{Dielectric universality}    \label{sec_dielectric_universality}
        \input{universal_dielectric_response/sec_dielectric_universality}
    \subsection{Universal dielectric response as a fractional process}  \label{sec_DU_as_fractional}
        \input{universal_dielectric_response/sec_DU_as_fractional}
        
\section{Numerical approximation of fractional integrals}   \label{sec_numerical_approximation_FI}
    \subsection{Fundamentals}   \label{sec_fundamentals}
        \input{numerical_approximation/sec_fundamentals}
    \subsection{Linear multi-step method}   \label{sec_LMSM}
        \input{numerical_approximation/sec_LMSM}
    \subsection{Classical convolution quadratures}  \label{sec_classical_conv_quad}
        \input{numerical_approximation/sec_classical_conv_quad}
    \subsection{Fractional quadratures and the Gr\"unwald-Letniknov fractional derivative} \label{sec_frac_quad_GL_derivative}
        \input{numerical_approximation/sec_frac_quad_GL_derivative}
    \subsection{Gamma and factorial functions overflow and  the short-memory principle for Gr\"unwald--Letnikov approximation}   \label{sec_gamma_overflow}
        \input{numerical_approximation/sec_gamma_overflow}
    \subsection{Interpolation methods and the Fractional Newton-Cotes formula} \label{sec_interpolation_methods}
        \input{numerical_approximation/sec_interpolation_methods}
    \subsection{Numerical examples} \label{sec_numerical_examples}
        \input{numerical_approximation/sec_numerical_examples}
        
\section{Conclusions}   \label{sec_conclusions}
    \input{conclusions/sec_conclusions}
    
\section{Acknowledgements}  \label{sec_acknowledgements}
    \input{acknowledgements/sec_acknowledgements}

\input{references/references}
\end{document}

%% file: abstract/sec_abstract.tex
The aim of this work is to introduce the main concepts of Fractional Calculus, followed by one of its applications to classical electrodynamics, illustrating how non-locality can be interpreted naturally in a fractional scenario. In particular, a result connecting fractional dynamics to high frequency dielectric response is used as motivation. In addition to the theoretical discussion, a comprehensive review of two numerical procedures for fractional integration is carried out, allowing one immediately to build numerical models applied to high frequency electromagnetics and correlated fields. 

%% file: introduction/sec_introduction.tex
Fractional calculus has a long history, at least as old as the usual calculus \cite{ross_1975}. Concepts on derivatives and integrals of general order have remained in the realm of pure mathematics for a long time, but have been applied to a variety of different fields in the last decades, ranging from control theory to electrodynamics,  economics and quantum mechanics \cite{hermann_2014}. 
Among several applications and interpretations, fractional calculus (FC) seems to constitute a well-suited apparatus to model systems with memory, i.e., systems whose response to a local stimulus depends on their whole ``history'' \cite{hermann_2014,tarasov_2013}. Such property makes FC an important tool in linear classical electrodynamics, since field analysis inside bulk matter takes into account non-local contributions of the stimulus, both in time and space \cite{oughstun_2006}. In fact, there is a formal similarity between fractional integral development and Maxwell equations for non-local media, since in both scenarios convolution integrals play the main role \cite{tarasov_2013, oughstun_2006, hilfer_2000}, connecting stimulus and response fields. In addition, there are evidences that for sufficiently high frequencies, dielectrics of completely distinct kinds present a power-law response. This so-called dielectric universal response was discovered first by Jonscher \cite{jonscher_1977, jonscher_1999} rises naturally in theory when one uses fractional integrals in Maxwell non-local equations, as shown by Tarasov in \cite{tarasov_2008, tarasov_2010, tarasov_2009}.

It is natural that, in face of the potential applicability of FC to classical electrodynamics, numerical methods are devised to approximate solutions of complex problems. More essentially, FC machinery is not trivial\footnote{Very often, integral or derivative of simple functions such as $\exp(t)$ lead to results expressed by unusual special functions. Vide, e.g., Hermann \cite{hermann_2014}.}, so that numerical methods might be useful even for simple estimates. In this context, the present work intends to present the fundamentals of FC, as well on numerical methods for fractional integral computation. Such concepts are useful for those who seek for alternative formulations of high-frequency response or, in a wide sense, to any non-local response theory \cite{tarasov_2013, tarasov_2010}.

The article is organized as follows:
\begin{itemize}
    \item   In Section \ref{sec_fractional_calculus} FC fundamentals are exposed, ranging from basic definitions up to identification criteria for fractional operators. The development focuses on two main definitions: Riemann-Liouville integral and Gr\"unwald-Letnikov derivative, due to their interconnection.  
    
    \item   Section \ref{sec_universal_dielectric_response} carries out a review on non-local fields in classical linear electrodynamics and dielectric response, where Jonscher's dielectric universality is discussed.

    \item   In Section \ref{sec_numerical_approximation_FI} a constructive review on numerical methods for fractional integration is carried out. The analysis focuses on two different tools: linear multi-step and Newton-Cotes methods, discussed in both classical and fractional scenarios. Section \ref{sec_numerical_examples} presents a few illustrations of the preceding methods.
\end{itemize}

A deep treatise on the mathematical aspects of fractional calculus is found in Samko et al.\cite{samko_1993}, whereas more applied approaches are given by Hermann \cite{hermann_2014}, Oldham and Spanier \cite{oldham_2006} and Tarasov \cite{tarasov_2010}. An introductory course on linear electrodynamics can be found in Greiner \cite{greiner_1998}, whereas a much more detailed development on field response can be found in Oughstun \cite{oughstun_2006}. Classical texts such  as Dahlquist \cite{dahlquist_2003}, Ralston and Rabinowitz \cite{ralston_2001} and Golub and Ortega \cite{golub_1991} provide the whole background necessary to the understanding of the considered numerical methods, at least in the classical scenario. The generalization for fractional integrals through linear multi-step methods is strongly based on \cite{lubich_1986}, supported by the formalism introduced in \cite{matthys_1976, wolkenfelt_1981, lubich_1988a,lubich_1988b}. On the other hand, fractional Newton-Cotes formulae development can be found, for example, in Li and Zeng \cite{li_2015}. 

%% file: fractional_calculus/sec_fundamental_concepts.tex
Usual differentiation and integration are well-known concepts in science and technology, and have a strong theoretical basis, built on the efforts of Newton, Leibniz, Weirstrass, among many other. However, as old as the usual differential and integral calculus, there is the so-called fractional calculus apparatus. According to Ross \cite{ross_1975}, fractional calculus origins goes back to a dialog between L’Hospital and Leibniz, in which the meaning of the operation $d^{1/2}/dt^{1/2}$  was discussed. FC main idea consists in generalizing the order of a derivative,say $n$, allowing one to perform derivations -- and also integrations -- of $\al$-order, with $\al \in \conjC$. Immediately, one may ask about the meaning of the operation
\begin{equation}\label{eq_fractional_derivative}
    \frac{d^{\al}}{dt^{\al}}f(t), ~ \al \in \mathbb{C}
\end{equation}
and how it is operationalized. An intuitive definition of fractional derivative was given by Fourier \cite{hermann_2014, hilfer_2000}: it is well known that Fourier transform maps the derivative of its operand, say $f(t)$, to a product between Fourier transform of the original function and a power of $\jmath \w$, with $\jmath = \sqrt{-1}$, i.e., 
\begin{equation}\label{eq_FT_of_derivative}
    \ves{F}\left[   \frac{d^{n}}{dt^{n}}f(t)     \right](\w) = (\jmath \w)^{n} \ves{F}[f(t)], \quad \text{with } n \in \conjN.
\end{equation}
Generalizing the derivative order to $\al \in \conjC$ and performing an inversion, one obtains
\begin{equation}\label{eq_fourier_FD}
    D^{\al}[f](t) =  \frac{d^{\al}}{dt^{\al}}f(t) = \frac{1}{2 \pi}\int\limits_{-\infty}^{\infty}{(\jmath \w)^{\al}\ves{F}(\w) d\w}, 
\end{equation}
where $\ves{F}(\w) = \int{f(t)\exp(-\jmath \w t)dt}$ is the Fourier transform of $f$.
For example, if $f(t) = \sin(\w_0 t)$, one obtains $D^{\al}[f](t) = |\w_0|^{\al}\sin\left(\w_0 t + \frac{\pi}{2}\al\right)$, which leads to the usual result when $\al \in \conjZ$. This is not a coincidence, but a fundamental requirement for fractional derivatives, as will be discussed later.

Despite its simplicity, Fourier definition of fractional derivative is not unique: there are many others. For instance, consulting \cite{hermann_2014}, \cite{hilfer_2000} and \cite{oliveira_2014}, one finds more than thirty  definitions, not generally equivalent. Before introducing two of them,  it is important to emphasize the fact that Fourier fractional derivative, Eq. \eqref{eq_fourier_FD}, is an integral in the classical sense. This is a general characteristic of fractional operators: derivatives and integrals are deeply connected, similarly to what occurs in the usual scenario, where Fundamental Theorem of Calculus relates both operations. To visualize such connection in the fractional scenario,  we first remind that an arbitrary number $\al$ can be decomposed on its integer and fractional parts, i.e., $\al = [\al] + \{\al\}$, with $[\al] \in \conjZ$ and $\al \in [0, 1[$. Therefore,
\begin{equation}\label{eq_semigroup_FD}
    D^{\al}[f](t) = D^{ [a] + \{\al\}}[f](t) = 
    \begin{cases}
         D^{[\al]}D^{\{\al\}}[f](t); ~ \text{or}  \\
         \\
         D^{\{\al\}}D^{[\al]}[f](t).
    \end{cases}
\end{equation}
Eq. \eqref{eq_semigroup_FD} uses the important semi-group property, which basically ensures the commutativity of $D[\cdot]$. In order to connect Eq. \eqref{eq_semigroup_FD} to a fractional integral \cite{hermann_2014}, one expands $D^{\al}$ into $D^{n}D^{\al-n}$, with  $n = [\al] + 1$:
\begin{eqnarray}\label{eq_connection_FD_FI}
    D^{\al}[f](t)   &   =   &   D^{\al - n + n}[f](t)       \nonumber   \\
                    &   =   &   D^{n}D^{\al-n}[f](t)        \nonumber   \\
                    &   =   &   D^{n}D^{\{\al\} - 1}[f](t)  \nonumber   \\
                    &   =   &   D^{n}I^{1-\{\al\}}[f](t), 
\end{eqnarray}
that is, an $\al$-order fractional derivative is related to an  $(1-\al)$-order iterated fractional integral, denoted by $I[\cdot]$, followed by an integer n-order derivative, with $n = [\al] + 1$. This intricate connection has led Oldham and Spanier \cite{oldham_2006} to coin picturesque terms such ``differintegral" and ``differentegration", since if one defines an operation identified as a fractional integral, its derivative is readily provided, according Eq.\eqref{eq_connection_FD_FI}. In fact, for most of the definitions of fractional operators, such as the Caputo and Riemann-Liouville approaches, one starts with the fractional integral, since the derivative counterpart is already defined. One remarkable exception is the Gr\"unwald-Letnikov derivative, which is based on the generalization of finite differences.

In what follows we shall focus on two definitions: Riemann-Liouville and Gr\"unwald-Letnikov. These two approaches were chosen due to their interesting interconnection, of great importance for numerical purposes.

%% file: fractional_calculus/sec_RL_approach.tex
The following result, due to Cauchy\footnote{According Ross \cite{ross_1975}, the result shown in Eq. \eqref{eq_iterated_integral} is due to Dirichlet.} \cite{hilfer_2000, hermann_2014, samko_1993} relates the iterated integral
\begin{equation}\label{eq_iterated_integral}
    {}_{a}I^{n}[f](t) = \int\limits_{a}^{t}{dt_{n-1}    \int\limits_{a}^{t_{n - 1}}{dt_{n_2}{    \hdots     \int\limits_{a}^{t_2}{dt_1  \int\limits_{a}^{t_1}{dt_0 f(t_0)}    }  }}    }
\end{equation}
to a simpler form, given by
\begin{equation}\label{eq_cauchy_iterated_integral}
    {}_{a}I^{n}[f](t) = \frac{1}{(n-1)!}\int\limits_{a}^{t}{(t-t')^{n-1}f(t')dt'}.
\end{equation}
Generalizing the order in Eq. \eqref{eq_cauchy_iterated_integral} to an arbitrary $\al$, one obtains left and right-sided Riemman-Liouville integrals of order $\al > 0$, respectively: 
\begin{equation}\label{eq_LSRLI}
    {}_{a}I^{\al}_{+}[f](t) = \frac{1}{\G(\al)}\int\limits_{a}^{t}{\frac{f(t')}{(t-t')^{1-\al}}dt'}, ~ \text{with } t > a,
\end{equation}
and
\begin{equation}\label{eq_RSRLI}
    {}_{a}I^{\al}_{-}[f](t) = \frac{1}{\G(\al)}\int\limits_{t}^{a}{\frac{f(t')}{(t'-t)^{1-\al}}dt'}, ~ \text{with } t < a, ~ \text{and}.
\end{equation}
In both Eqs. \eqref{eq_LSRLI} and \eqref{eq_RSRLI},  $\G(\cdot)$ refers to the Gamma function. According to Eq. \eqref{eq_connection_FD_FI}, to obtain the Riemann-Liouville fractional derivatives, it suffices to perform an integer differentiation of order  on the expressions above.

%% file: fractional_calculus/sec_GL_approach.tex
It is well-known that the usual $n$-order derivative of  $f(t)$ can be written as the limit
\begin{equation}\label{eq_classical_deriv_FD}
    \lim\limits_{h \to 0}{\frac{\Del^{n}_{h}[f](t)}{h^n}} = \lim\limits_{h \to 0}{ \frac{\sum\limits_{k = 0}^{n}{(-1)^{k}\binom{n}{k}}f(t-kh)}{h^n}  },
\end{equation}
in which
\begin{equation}\label{eq_binom_coeff}
    \binom{n}{k} = \frac{n!}{k!(n-k)!}
\end{equation}
is the binomial coefficient. Suggestively, one might generalize the difference operator $\Del_{h}^{n}[\cdot]$ to
\begin{equation}\label{eq_general_FD_operator}
    \Del_{h}^{\al}[f](t) = \sum\limits_{k = 0}^{\infty}{\binom{\al}{k}f(t-kh)}, ~ \text{with } \al > 0.  
\end{equation}
Now, the binomial coefficient is given by
\begin{equation}\label{eq_general_binom_coeff}
    \binom{\al}{k} = \frac{\G(\al+1)}{\G(k+1)\G(\al-k+1)}.
\end{equation}
The Gr\"unwald-Letnikov fractional derivative is given by the following limit:
\begin{equation}\label{eq_GL_frac_deriv}
    D^{\al}[f](t) = \lim\limits_{h \to 0}{\frac{\sum\limits_{k = 0}^{\infty}{\binom{\al}{k}f(t-kh)}}{h^{\al}}}.
\end{equation}
If $h > 0$, one deals with backward differences, otherwise one has forward difference. Care must be taken when $\al < 0$, once the series given by Eq. \eqref{eq_general_FD_operator} can easily diverge. However, as pointed out by Samko \emph{et al.}\cite{samko_1993}, the series shown in Eq. \eqref{eq_general_FD_operator} eventually converges if $f(t)$ is a non-periodic function which smoothly -- and fastly -- decays to zero at infinity. Within these conditions, Eq. \eqref{eq_GL_frac_deriv} turns into Gr\"unwald-Letnikov fractional integral \cite{samko_1993, tarasov_2013}.

%% file: fractional_calculus/sec_identification_criteria.tex
As pointed out earlier, the number of different definitions for fractional integrals and derivatives is remarkable. Surprisingly, different definitions lead to different results when applied on the same operand. In light of such plurality, it is natural to ask how can one identify a ``legitimate" fractional operator. In 1975, Ross \cite{ross_1975} proposed five criteria in order to ensure the fractionality of a given operator:
\begin{enumerate}
    \item   The $\al$-order derivative of an analytic function $f(z)$  must be analytic on $z$ and $\al$;
    \item   If $\al \in \conjZ^{+}$, the result must be the same obtained by usual -- integer order -- differentiation. If $\al \in \conjZ^{-}$, the result must be the same of the $\al$-fold iterated integral;
    \item   The operator must be linear;
    \item   If $\al = 0$, the operand must remain unaltered. In other words, 0 is the neutral element;
    \item   Given $\al, \be \in \conjC$, the semi-group property must hold:
            \begin{equation}\label{eq_general_semigroup_prop}
                D^{\al}D^{\be}[f] = D^{\al+\be}[f];
            \end{equation}
\end{enumerate}
In addition to the criteria proposed by Ross, there is a proof given by Tarasov \cite{tarasov_2013a} which shows that for a legitimate fractional operator, the usual Leibniz rule, given by
\begin{equation}\label{eq_usual_leibniz_rule}
    \frac{d^n}{dt^n}[fg](t) = \sum\limits_{k = 0}^{n}{\binom{n}{k}f^{(n-k)}(t)g^{(k)}(t)}, ~\text{with } n \in \conjZ
\end{equation}
\emph{must} be broken.
Ortigueira and Machado \cite{ortigueira_2015} used this result to improve Ross criteria, stating that
\begin{enumerate}   \setcounter{enumi}{5}
    \item   For a true fractional operator, the generalized Leibniz product rule must hold:
            \begin{equation}\label{eq_fractional_leibniz_rule}
                D^{\al}[fg](t) = \sum\limits_{k = 0}^{\infty}{\binom{\al}{k}D^{\al-k}[f](t) D^{k}[g](t)}, \al \in \conjR.
            \end{equation}
\end{enumerate}
In this same work, the authors show that Riemann-Liouville and Gr\"unwald-Letnikov fractional derivatives – among several other – obey Eq. \eqref{eq_fractional_leibniz_rule}.
The six criteria shown above work as a sort of filter, since there is a profusion of allegedly fractional operators in the literature.

%% file: universal_dielectric_response/sec_linear_media.tex
In an arbitrary material media, Maxwell equations are given by
\begin{equation}\label{eq_maxwell}
    \begin{cases}
        \diver{\varvet{D}}(\rr,t)   & = \quad    \rho(\rr,t) \nonumber   \\
        \rot{\varvet{E}}(\rr,t)     & = \quad   -\partial_t \varvet{B}(\rr,t)   \nonumber   \\
        \rot{\varvet{H}}(\rr,t)     & = \quad   \varvet{j}(\rr,t) + \partial_t \varvet{D}(\rr,t)    \nonumber   \\
        \diver{\varvet{B}}(\rr,t)   & = \quad   0,
    \end{cases}
\end{equation}
in which $\rho$ and $\varvet{j}$ are the free charge and current densities, $\varvet{E}$, $\varvet{D}$, $\varvet{B}$ and $\varvet{H}$ are the electric field, electric displacement, magnetic induction and magnetic field, respectively.  The fields  $E$, $P$ and $D$ are related by
\begin{equation}\label{eq_constitutive}
    \varvet{D}(\rr,t) = \eps_0 \varvet{E}(\rr,t) + \varvet{P}(\rr,t), 
\end{equation}
in which  $\varvet{P}$ is the polarization of the medium and $\eps_0$ is the vacuum permittivity.  The simplest interpretation for $\varvet{P}$ states that this field represents the dipole moment per unit (of a non-infinitesimal) volume. For linear media, there is a direct relation between $\varvet{E}$ and $\varvet{P}$, expressed as
\begin{equation}\label{eq_constitutive_linear}
    \varvet{P}(\rr,t) = \eps_0 \chi \varvet{E}(\rr, t)
\end{equation}
in which  $\chi$ is the electric susceptibility. Combining Eqs. \eqref{eq_constitutive} and \eqref{eq_constitutive_linear}, one obtains the so-called constitutive relation for linear media:
\begin{equation}\label{eq_constitutive_linear_final}
    \varvet{D}(\rr,t) = \eps \varvet{E}(\rr,t),
\end{equation}
where the permittivity is given by $\eps \equiv \eps_0 (1 + \chi)$. One should notice that we are assuming isotropic media, since $\eps$ is a scalar function\footnote{In the most general case, $\eps$ is a tensorial function.}.

The simplicity of Eq. \eqref{eq_constitutive_linear_final} is due to two strong implicit assumptions:
\begin{itemize}
    \item   a stimulus provided by $\varvet{E}$ at a given spatial coordinate $\rr'$ affects $\varvet{D}$ only at $\rr'$; and
    \item   a stimulus provided by $\varvet{E}$ at a given time instant $t'$ is immediately felt by $\varvet{D}$.
\end{itemize}
In fact, in the general case, linearity between $\varvet{D}$ and $\varvet{E}$ is expressed by \cite{oughstun_2006}, \cite{greiner_1998}
\begin{equation}\label{eq_general_constitutive}
    \varvet{D}(\rr,t) = \int\limits_{-\infty}^{\infty}{  \int\limits_{-\infty}^{\infty}{\eps(\rr',t'; \rr, t)}\varvet{E}(\rr',t')d^{3}\rr'dt'},
\end{equation}
where $\eps(\cdot)$ is the integral kernel. If the medium is local \cite{oughstun_2006}, assumption 1 above is valid, since
\begin{equation}\label{eq_local_eps}
    \eps(\rr',t';\rr,t) = \eps(t',t)\del(\rr-\rr')
\end{equation}
and then
\begin{equation}\label{eq_local_D}
    \varvet{D}(\rr,t) = \int\limits_{-\infty}^{\infty}{\eps(\rr;t',t)\varvet{E}(\rr,t')dt'}.
\end{equation}
Spatial locality roughly means that the medium molecular constituents are independent from each other \cite{oughstun_2006}. Such assumption is, in some measure, admissible. On the other hand, assumption 2 is not realistic: it assumes that information propagates at infinity speed in the medium. In other words, it violates the principle of causality. However, in static and quasi-static scenarios, Eq. \eqref{eq_constitutive_linear_final} is a good approximation \cite{greiner_1998}.

If we make a further assumption, assuming that the medium is temporally homogeneous\footnote{Temporal homogeneity is a quite reasonable assumption, since it means that the functional structure of $\eps$ doesn’t change in time. In other words, the system is time-shift invariant.}, we have
\begin{equation}\label{eq_LSI_eps}
    \eps(\rr; t, t') =  \eps(\rr; t - t'),
\end{equation}
which leads to a temporal convolution between $\varvet{D}$ and $\varvet{E}$:
\begin{equation}\label{eq_conv_eps_E}
    \varvet{D}(\rr,t) = \int\limits_{-\infty}^{\infty}{\eps(\rr,t-t')E(\rr,t')dt' \equiv \left( \eps(\rr) * \varvet{E}(\rr)     \right)(t)}.
\end{equation}
Eq. \eqref{eq_conv_eps_E} makes explicit the fundamental characteristic of temporal non-locality, since the response $\varvet{D}$ depends, in a certain manner dictated by $\eps$, on the whole history of the stimulus $\varvet{E}$. A last improvement of Eq. \eqref{eq_conv_eps_E} requires us to take into account the principle of causality: since response can not precede  stimulus, we must have  $\eps(\rr;t-t') = 0$ when $t-t' < 0$. Therefore, one deals now with a causal temporal convolution:
\begin{equation}\label{eq_causal_conv_eps_E}
%     \varvet{D}(\rr,t) = \int\limits_{-\infty}^{t}{\eps(\rr,t-t')}
    \varvet{D}(\rr,t) = \int\limits_{-\infty}^{t}{\eps(\rr,t-t')E(\rr,t')dt'}.
\end{equation}
A medium in which Eq. \eqref{eq_causal_conv_eps_E} is valid is denominated as linear temporally dispersive \cite{oughstun_2006}.

Applying Fourier transform on both sides of Eq. \eqref{eq_causal_conv_eps_E} leads us to
\begin{equation}\label{eq_FT_causal_conv_eps_E}
    \varvet{D}(\varvet{k},\w) = \left[\eps(\varvet{k},\w) \varvet{E}(\varvet{k},\w)  \right] * \left(   \pi \del(\w) + \frac{\jmath}{\w}       \right),
\end{equation}
in which $\varvet{k}$ and $\w$ are the wave vector and angular frequency, respectively.  The convolution in Eq. \eqref{eq_FT_causal_conv_eps_E} is due to the fact that in order to ensure causality, we have performed a product in time-space between $\eps(\rr,t)\varvet{E}(\rr,t)$ and the Heaviside function, $\theta(\cdot)$,  whose Fourier transform pair is given by
\begin{equation}\label{eq_FT_step_function}
    \theta(t) = 
    \begin{cases}
        1, & t \geq 0 \\
        0, & t <0
    \end{cases}
    \quad \Leftrightarrow \quad
    \theta(\w) = \pi \del(\w) + \frac{\jmath}{\w}.
\end{equation}
Under adequate assumptions of symmetry and analyticity  of stimulus and response functions, one obtains from Eq. \eqref{eq_FT_causal_conv_eps_E} the important Kramers-Kronig (KK) relations \cite{oughstun_2006, greiner_1998}, which connect real and imaginary parts of the system’s optical response. KK relations are built on the theory of Hilbert transformation, and the interested reader may consult  the works of Oughstun \cite{oughstun_2006}, Greiner \cite{greiner_1998} and Witthaker et al. \cite{whittaker_2015} for further details.

%% file: universal_dielectric_response/sec_empirical_models.tex
Once the role played by temporal dispersion was made explicit, one might ask about the behavior of the field $\varvet{P}$  when subjected to a stimulus from $\varvet{E}$ . There are many empirical models which relate these fields \cite{oughstun_2006, garrappa_2016}. In general lines, one tries to find a phenomenological model for an individual dipole moment $\varvet{p}$ and then computes its average structure, related to the macroscopic polarization $\varvet{P}$:
\begin{equation}\label{eq_macro_P}
    \varvet{P}(\rr, t) = \langle \varvet{p}(\rr,t) \rangle = \sum \limits_{i}{f_i \varvet{p_i}(\rr,t)}, 
\end{equation}
in which the brackets indicate ensemble average and  $f_i$ are dipolar statistical weights. For example, a working model for non-dense media \cite{greiner_1998} starts from the hypothesis that each electron of the medium is subject to classical damping and elastic forces. Therefore, one has the movement equation
\begin{equation}\label{eq_damped_oscillator}
    \ddot{\varvet{x}}_i + \g_i \dot{\varvet{x}}_i + \w_i \varvet{x}_i = \frac{e}{m}\varvet{E},
\end{equation}
in which $\g_i$, $\w_i$ and $m$ are the i-th damping constant, characteristic frequency and electronic mass, respectively. The i-th state is characterized by $\g_i$ and $\w_i$, i.e., $f_i$ equals the number of electrons with damping constant and angular frequency $\g_i$ and $\w_i$. If the field is harmonic, i.e., $\varvet{E}(t) = \varvet{E_0}\exp(-\jmath \w t)$, Eq. \eqref{eq_damped_oscillator} has a solution given by
\begin{equation}\label{eq_solution_damped_oscillator}
    \varvet{p}_i (t) = \frac{e^2}{m}\frac{1}{({\w_i}^2 - \w^2) - \jmath \g_i \w}\varvet{E}(t).
\end{equation}
Using Eq. \eqref{eq_solution_damped_oscillator} in Eq. \eqref{eq_macro_P} and reminding that $\varvet{p} = e \varvet{x}$, one obtains
\begin{equation}\label{eq_macro_P_damped_oscillator}
    \varvet{P}(\w) = \frac{N e^2}{m}\sum\limits_{i}{\frac{f_i}{({\w_i}^2 - \w^2) - \jmath \g_i \w}\varvet{E}(\w)},
\end{equation}
where $N$ is the number of molecules per unit volume. With Eq. \eqref{eq_macro_P_damped_oscillator} and \eqref{eq_FT_causal_conv_eps_E}, one finds expression for the optical responses  $\eps(\w)$ or $\chi(\w)$ .

%% file: universal_dielectric_response/sec_dielectric_universality.tex
As pointed out at the beginning of the section, there are several models for polarization, such as the Clausius-Massoti and Cole-Cole models. Interested reader should refer to \cite{oughstun_2006} and \cite{garrappa_2016} for deeper analysis on the subject. However, it is interesting to mention the important model due to Debye \cite{oughstun_2006, jonscher_1999}. In this model, the electrical dipole $\varvet{p}$ is supposed to suffer viscous resistance, being allowed to rotate in the medium, but it is not able to interact to other dipoles and there is not Brownian motion in such ``fluid'' medium. This model leads to the following movement equation [4]:
\begin{equation}\label{eq_debye_model}
    \dot{\varvet{p}} + \frac{1}{\tau}\varvet{p} = a \varvet{E},
\end{equation}
in which $a$ is a coupling constant and $\tau$ is the molecular relaxation time, i.e., the characteristic time interval the dipole takes to return to a random orientation after $\varvet{E}$ ceases. It can be shown \cite{oughstun_2006, jonscher_1977} that Debye model leads to a susceptibility of the form
\begin{equation}\label{eq_electric_susc_debye}
    \chi(\w) = \frac{1}{\eps_0}\frac{N a \tau}{1 - \jmath \w \tau}.
\end{equation}
It is important to emphasize that Debye’s model can be improved \cite{oughstun_2006, garrappa_2016}, but still fails at some point. Experimental observations in this direction can be found in  \cite{zhao_2016} and \cite{gao_2015}, for instance. In fact, one might build more complete models for dielectric response, but at the expense of dealing with more variables and a increasing degree of complexity. Instead,  we look now at the behavior of a large class of media at high frequencies. Interestingly, very distinct materials show some power law behavior for the dielectric response. Let us denote the susceptibility by
\begin{equation}\label{eq_complex_susc}
    \chi(\w) = \chi'(\w) + \jmath \chi''(\w).
\end{equation}
According Jonscher \cite{jonscher_1977, jonscher_1999}, the following proportionality is observed at high frequencies:
\begin{equation}\label{eq_universal_susc}
    \chi(\w) \propto (\jmath \w)^{n - 1} ~~ \text{with } n \in ]0,1[ ~\text{and } \w \gg \frac{1}{\tau}.
\end{equation}
Eq. \eqref{eq_universal_susc} shows that the ratio between imaginary and real parts of $\chi$ does not depend on the frequency:
\begin{equation}\label{eq_complex_susc_ratio}
    \frac{\chi''(\w)}{\chi'(\w)} = \cot\left(n \frac{\pi}{2} \right).
\end{equation}
In \cite{jonscher_1999, tarasov_2008, tarasov_2010} some theoretical scenarios are discussed in order to explain this universal behavior. Some of them are built on generalizations of the Debye model; for instance, models where there are many relaxation times, instead of only one. In such models, one should consider not only one relaxation time, but an ensemble average on the distribution. In addition, the concept of fractionality is also considered, i.e., according to such approach, the power-law behavior could be associated to some sort of scale invariance of the medium. 

%% file: universal_dielectric_response/sec_DU_as_fractional.tex
Among a plurality of interpretative scenarios for universal dielectric response, we turn our focus to the formal aspects subjacent to the phenomenon. Following Tarasov in \cite{tarasov_2008, tarasov_2010, tarasov_2009}, we start from constitutive relations for linear media and the high-frequency universal behavior, Eq. \eqref{eq_universal_susc}, moving towards a fractional integral relation.

According to the discussion carried out in Section \ref{sec_linear_media}, one has to expect the following relation between $\varvet{P}$ and $\varvet{E}$:
\begin{equation}\label{eq_constitutive_P_E}
    \varvet{P}(\rr,t) = \eps_0 \int\limits_{-\infty}^{\infty}{\chi(t-t')\varvet{E}(\rr,t')dt'}.
\end{equation}
Eq. \eqref{eq_constitutive_P_E} assumes that the susceptibility $\chi$ does not depend on $\rr$ and is itself a causal function, so that the integral limits can be set at infinity. Taking temporal Fourier transform on both sides, one has
\begin{equation}\label{eq_FT_constitutive_P_E}
    \varvet{P}(\rr,\w) = \eps_0 \chi(\w)\varvet{E}(\rr,\w).
\end{equation}
Using Eq. \eqref{eq_universal_susc} as an equality with the proportionality constant equals one in the previous result, and denoting $-\al = n - 1$, we obtain
\begin{equation}\label{eq_macro_P_E}
    \varvet{P}(\rr,\w) = \eps_0(\jmath \w)^{-\al}\varvet{E}(\rr,\w), ~ \text{with } \al \in ]0,1[ ~\text{and } \w \gg \frac{1}{\tau}.
\end{equation}
Samko \emph{et al.}\cite{samko_1993} show that the Fourier transform derivative property, given by Eq. \eqref{eq_FT_of_derivative} remains valid in the fractional scenario. The right-hand side of Eq \eqref{eq_macro_P_E} contains a term $(\jmath \w)^{-\al}$, leading to an $\al$-order integral in the time domain \cite{tarasov_2009}:
\begin{equation}\label{eq_fractional_P_E}
    \varvet{P}(\rr,t) = \eps_0 I^{\al}_{+}[\varvet{E}](\rr,t) ~\text{with } \al \in ]0,1[.
\end{equation}
The result shown in Eq. \eqref{eq_fractional_P_E} allows an interesting interpretation: the polarization of the medium depends on $\varvet{E}$  past history in a more complex way than in the usual scenario. Temporal dispersion is related to memory effects, vide Eq. \eqref{eq_conv_eps_E}. But the fractional relation between stimulus and effect opens door to further interpretations; for instance, to those related to fractality and multiple time scales for the relaxation times.

One last comment: the previous development is due to Tarasov \cite{tarasov_2010, tarasov_2008, tarasov_2009}, but he was not the first to notice fractional characteristics in the universal dielectric response. In fact, in \cite{garrappa_2016}, Garrappa \emph{et al.} investigate several models whose original formulations made use of fractional integral kernels.

%% file: numerical_approximation/sec_fundamentals.tex
Once FC main ideas were exposed and one of its applications to electromagnetism was detailed, it is interesting to consider how one can discretize fractional integrals and derivatives and approximate them numerically. Obviously, such study would comprise an entire work by itself,  so that we will focus on the elements which provide the minimum apparatus on how to proceed with applied problems in FC using Riemann-Liouville integrals. The interested reader can refer to \cite{li_2015} and \cite{lubich_1986} for much deeper analysis.

In order to approximate the Riemann-Liouville (hereafter identified by RL) integral, one might explore its ``convolutional'' structure, shown in Eqs. \eqref{eq_LSRLI} and \eqref{eq_RSRLI}, since one has -- for the left-sided case --
\begin{equation}\label{eq_convolution_LSRLI}
    \frac{1}{\G(\al)}\frac{1}{t^{1 - \al}}*f(t),
\end{equation}
where * is the convolution symbol. An analog expression exists for the right-hand sided RL integral.

It may be tempting to approximate the convolution shown in Eq. \eqref{eq_convolution_LSRLI} in a straight manner, i.e, through discretization of the independent variable and then performing a discrete convolution in time domain or an usual product in the reciprocal space, through use of the discrete Fourier transform -- DFT. However, one should notice that the convolution kernel is singular at the origin when $\Re(\al) < 1$. Therefore, care must be taken in order to perform such approximation. 

%% file: numerical_approximation/sec_LMSM.tex
In this work we want to approximate left-sided RL integrals, Eq. \eqref{eq_LSRLI} with $a = 0$, through the use of fractional linear multistep method (FLMM), first proposed by Lubich \cite{lubich_1986}. This method is a generalization of the the well-known  class of linear multistep methods (LMM) used to approximate solutions of ordinary differential equations (ODEs) \cite{golub_1991}. Below, the main concepts of usual LMM are reviewed, building the path to the exposition of Lubich generalization.

Consider first the classical case in which $\al = 1$. One may employ LMMs to approximate the following  initial value problem:
\begin{equation}\label{eq_classical_IVP}
    \begin{cases}
        y'(t)    = f(t,y(t))      \\
        \\
        y(0)     = 0
    \end{cases}
    \Rightarrow \quad
    y(t) = \int\limits_{0}^{t}{f(t',y(t'))dt'}
\end{equation}
Let  $t_i = i \Del t$ be the i-th node of the mesh generated by discretization of the interval of interest and $y_i$ the approximate solution of Eq. \eqref{eq_classical_IVP} at $t_i$. For simplicity, we consider that the increment $\Del t$ is constant along the process. The simplest approach to approximate the considered ODE is the explicit Euler method \cite{golub_1991, ralston_2001}, given by
\begin{equation}\label{eq_explicit_euler}
    \begin{cases}
        y_{i+1} = y_i + \Del t f(t_i,y_i)   \\
        \\
        y_0 = \hat{y}_0
    \end{cases}
    \text{with} ~i \in [0, N-1].
\end{equation}
It can be shown \cite{golub_1991, ralston_2001} that explicit Euler scheme is a first order method, denoted by $\ves{O}(\Del t)$. However, the main aspect to be considered here is the fact that Eq. \eqref{eq_explicit_euler} is a one-step method, i.e., the value of $y_{i+1}$ depends only on the value of $y_i$. Another example of one-step method is the class of Runge-Kutta schemes. Alternatively, one might make use of other entries of $y$ to approximate the desired solution. In general, one would deal with the following difference equation
\begin{equation}\label{eq_LMSM}
    \sum\limits_{k = 0}^{N}{\al_k y_{n-k}} = \Del t\sum\limits_{k = 0}^{N}{\be_k f_{n-k}}
\end{equation}
where  $f_k = y'_k$ and $\{\al_k\}_{0}^{N}$ and $\{\be_k\}_{0}^{N}$  are two sets of constants that depend on the structure of the employed method. Eq. \eqref{eq_LMSM} is called a linear multi-step method. There is a plethora of methods based on it, including the Adams-Bashfort (explicit) and Adams-Moulton (implicit) methods. References \cite{golub_1991, ralston_2001} are recommended to the reader who seeks for detailed analysis on the subject.

It is clear that the coefficients $\al_k, \be_k$ play a fundamental role on both stability and convergence of the LMM shown in Eq. \eqref{eq_LMSM}. In fact, since one deals with a difference equation with constant coefficients, discrete Laplace transform rises as a powerful tool for analysis of such systems. It is defined by\footnote{One might notice the similarity between Eq. \eqref{eq_DLT} and the Z-transform of a causal signal, except from exponent’s signal. In fact, Z-transform is a useful resource in discrete system analysis and could be used to characterize the LMM given by Eq. \eqref{eq_LMSM} as well its counterpart given by Eq. \eqref{eq_DLT} \cite{coughlan_2007, galeone_2006}.}
\begin{equation}\label{eq_DLT}
    \ves{L}[y](z) = \sum\limits_{n = 0}^{\infty}{z^n y_n}
\end{equation}
and allows us to obtain the generating polynomials of the considered LMM \cite{ralston_2001, galeone_2006, lubich_1983, matthys_1976, wolkenfelt_1981}
\begin{equation}\label{eq_gen_poly}
    \begin{cases}
        \rho(z) &   = \sum\limits_{j = 0}^{N}{\al_j z^{n - j}}  \\
        \\
        \si(z)  &   = \sum\limits_{j = 0}^{N}{\be_j z^{n - j}}
    \end{cases}
\end{equation}
Since the polynomials $\rho(z)$ and $\si(z)$ completely determine the LMM, it is customary to designate the method by $(\rho, \si)$.

%% file: numerical_approximation/sec_classical_conv_quad.tex
As pointed out at the beginning of Section \ref{sec_fundamentals}, one might approximate the left-sided RL integral through the use of discrete convolutions. Despite the singularity at kernel's origin, Eq. \eqref{eq_convolution_LSRLI}, convolution is, in fact, a tool to be considered. We want to approximate an integral (possibly singular at some finite set) $I(t) = \int\limits_{0}^{t}{f(t')g(t-t')dt'}$ by means of a discrete convolution:
\begin{equation}\label{eq_discrete_conv}
    I(t = n\Del t) \equiv I_n \approx \sum\limits_{j = 0}^{n}{\w_j g_{n-j}} + \sum\limits_{j = 0}^{s}{\mu_{n j} g_{n-j}}, 
\end{equation}
where $\w_j$ depends on $\Del t$ and $s \leq n$ . The second sum is necessary to ensure convergence near the origin, and the coefficients $\mu_{n j}$ are the so-called starting quadrature terms \cite{lubich_1983, lubich_1988a, wolkenfelt_1981}. They can be obtained forcing the approximation given by Eq. \eqref{eq_discrete_conv} to be exact for polynomials of degree up to $s$ \cite{lubich_1988a}.

It can be shown \cite{lubich_1983, wolkenfelt_1981} that $\w_j$ is the j-th coefficient of the rational series
\begin{equation}\label{eq_rational_series}
    \w(z) = \frac{\si(z^{-1})}{\rho(z^{-1})} = \sum\limits_{j = 0}^{\infty}{\w_j z^j}.
\end{equation}
Observing Eq. \eqref{eq_rational_series}, one notices the need for a Laplace transform inversion in order to estimate $\w_j$. However, as pointed by Lubich in \cite{lubich_1988a} and \cite{lubich_1988b}, $\w_j$ is proportional to the s-order approximation of the inverse Laplace transform $\ves{L}^{-1}(s) = f(t)$, when $t$ is far from any singularity. In the general case, Eq. \eqref{eq_rational_series} requires a priori knowledge of the analytic Laplace transform of the considered LMM. If this is not the case, numerical procedures can be used to perform an inversion. 

%% file: numerical_approximation/sec_frac_quad_GL_derivative.tex
As considered in the previous section, linear multi-step methods can be used to approximate convolution integrals, in the classical sense. Interestingly, the generalization for the fractional case is quite straightforward: Lubich shows in \cite{lubich_1986} that the fractional left-sided Riemann-Liouville integral can be approximated through the coefficients of
\begin{equation}\label{eq_frac_coeffs}
    [\w(z)]^\al,
\end{equation}
in which $\w(z)$ is given by Eq. \eqref{eq_rational_series}. With this result, one has the approximation
\begin{equation}\label{eq_FLMSM}
    I^{\al}_{+}[f]_n \approx (\Del t)^{\al}\sum\limits_{j = 0}^{n}{\w^{\al}_{j}f_{n-j}} + (\Del t)^{\al}\sum\limits_{j = 0}^{s}{\mu_{nj}f_{n-j}}.
\end{equation}
One fundamental requirement about fractional linear multi-step methods (FLMM)  is that they must be implicit \cite{lubich_1986}.

A beautiful result rises when one uses the implicit Euler scheme, given by
\begin{equation}\label{eq_implicit_euler}
    y_i = y_{i - 1} + \Del t f_i.
\end{equation}
By direct inspection of Eq. \eqref{eq_implicit_euler} and application of Lubich fractional generalization, Eq. \eqref{eq_FLMSM}, one obtains the following polynomial:
\begin{equation}\label{eq_frac_binom}
    \w^{\al}(z) = (1 - z)^{-\al} = \sum\limits_{n = 0}^{\infty}{(-1)^{n}\binom{-\al}{n}z^n}
\end{equation}
Using this result in the discrete convolution given by Eq. \eqref{eq_FLMSM}, and discarding for the moment the starting weights, one has
\begin{equation}\label{eq_truncated_GL_deriv}
    I^{\al}_{+}[f]_n \approx (\Del t)^{\al}\sum\limits_{j = 0}^{n}{(-1)^j \binom{-\al}{j}f_{n-j}}
\end{equation}
which corresponds precisely to a truncated approximation for the Gr\"unwald-Letnikov derivative of $\al$-order, vide Eq. \eqref{eq_GL_frac_deriv}. From a theoretical perspective, Eq. \eqref{eq_truncated_GL_deriv} builds the path to the proof on the equivalence between Gr\"unwald-Letnikov (GL) and Riemann-Liouville fractional derivatives; see, for example, \cite{hilfer_2000} and \cite{lovoie_1976}. This connection allows us to approach theoretically a certain fractional problem using RL integrals, but numerically approximating it through GL derivatives of real negative orders \cite{li_2015}.

%% file: numerical_approximation/sec_gamma_overflow.tex
\begin{figure}[ht!]
    \centering
    \includegraphics[width = 0.7\columnwidth]{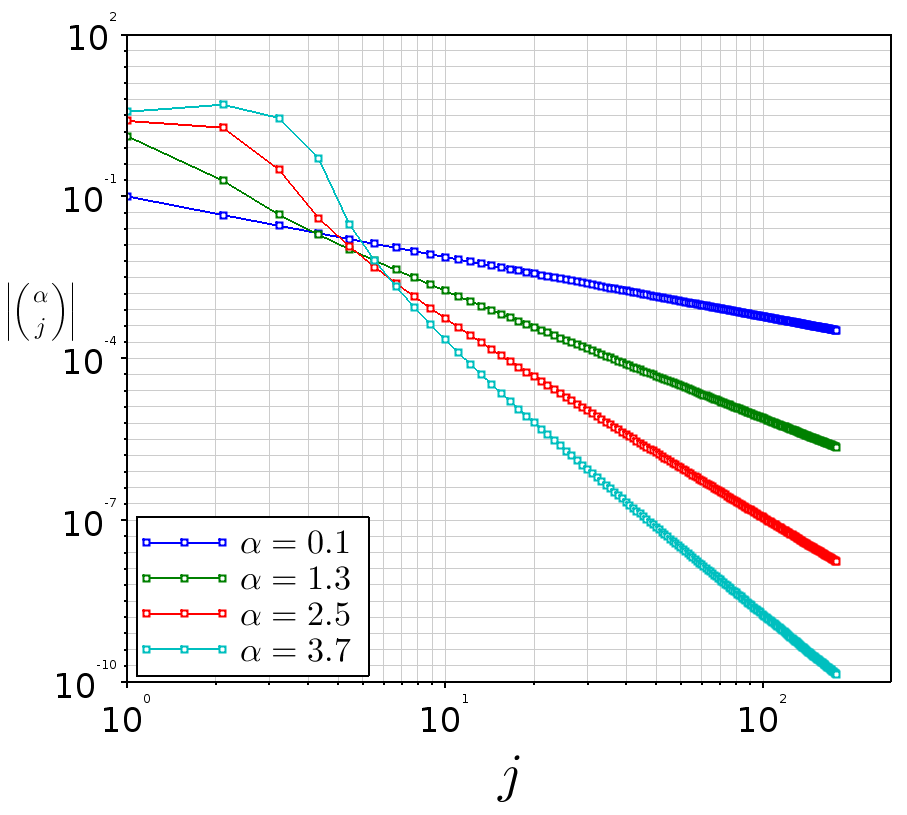}
    \caption{Evolution of the generalized binomial coefficient magnitude as $\al$ and $j$ increases.}
    \label{fig_binom_coeffs}
\end{figure}
\begin{figure}[ht!]
    \centering
    \includegraphics[width = 0.70\columnwidth]{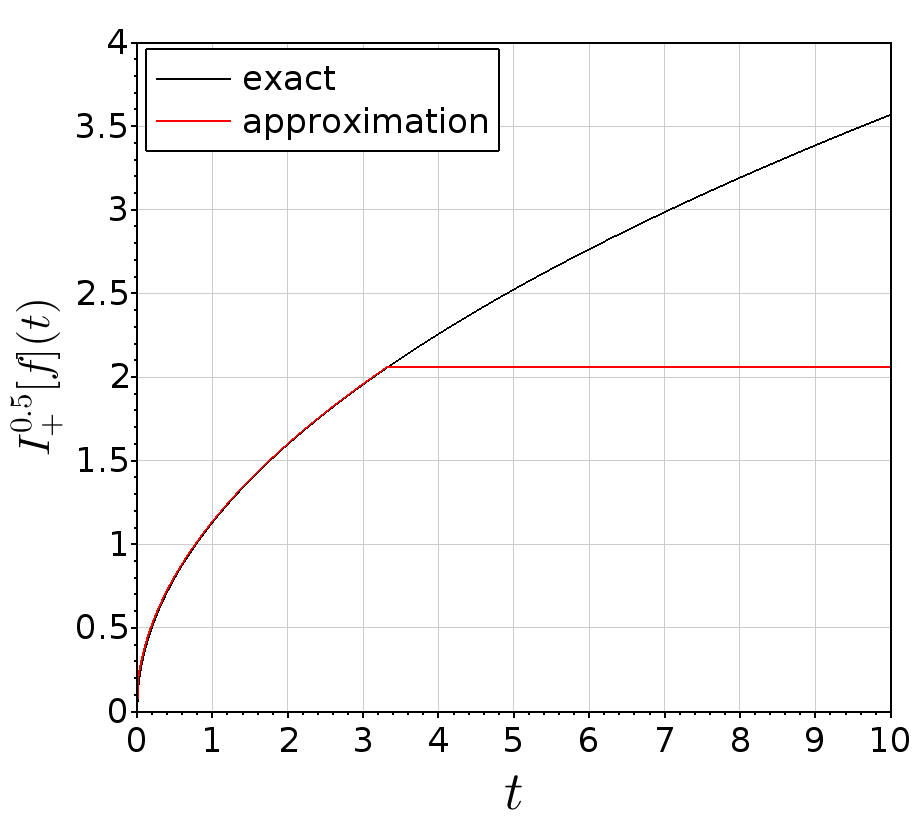}
    \caption{Artificial steadiness introduced in the Gr\"unwald-Letnikov approximation due overflow of gamma function for large $j$.}
    \label{fig_artificial_steadiness}
\end{figure}
\begin{figure}[ht!]
    \centering
    \includegraphics[width = 0.70\columnwidth]{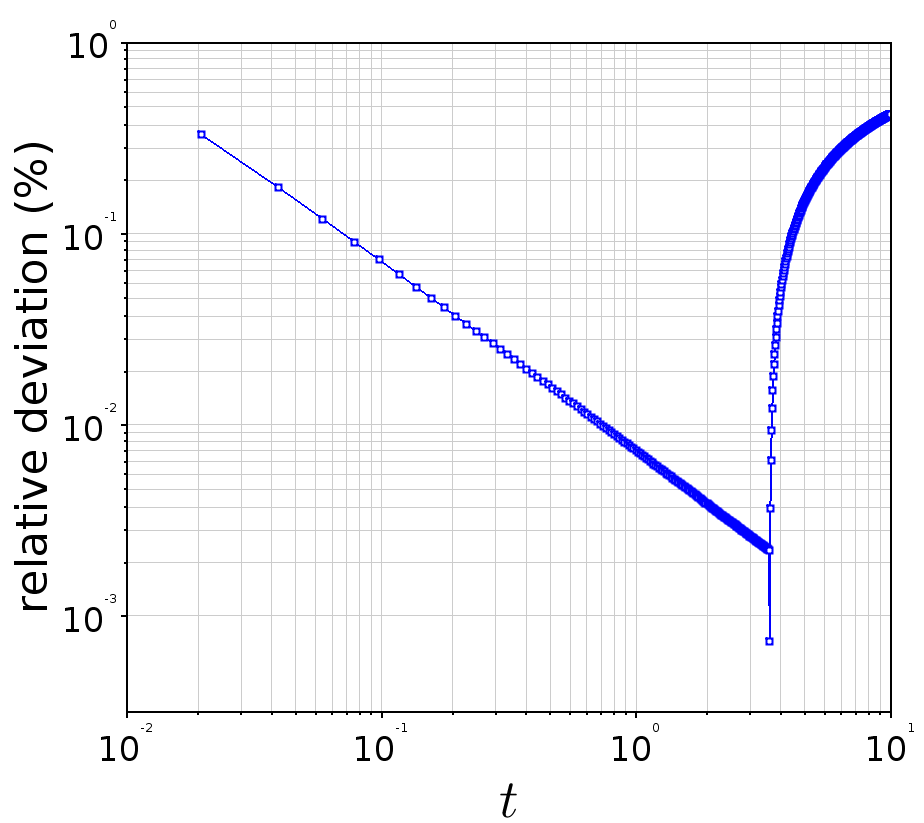}
    \caption{Error amplification in the Gr\"unwald-Letnikov approximation due to gamma function overflow.}
    \label{fig_error_amplification}
\end{figure}

The approximation obtained by means of Eq. \eqref{eq_truncated_GL_deriv} has an important drawback: the sequence 
\begin{equation}\label{eq_binom_coeffs_seq}
    \left\{\vert \w_j \vert =   \left \vert\binom{\al}{j}  \right \vert; \al > 0   \right\}
\end{equation}
decays rapidly to zero as $\al$ or $j$ increase, as can be seen in Figure \ref{fig_binom_coeffs}. Simultaneously, the coefficients $\binom{-\al}{j}$ diverge. In this sense, one loses any control over precision and stability of GL approximation.
Another error source due the form of Eq. \eqref{eq_binom_coeffs_seq} is related to the behavior of gamma and factorial functions for large numbers in computational environment. For instance, SCILAB language \cite{scilab_2012} returns infinity for $j! = \G(j+1); j > 170$. In practical terms, this overflow sets a limit on the mesh refinement. Such limitation can be visualized in Figure \ref{fig_artificial_steadiness}, where it is approximated the solution of $I^{0.5}_{+}[f](t)$ for $f(t) = 1, t \in ]0, 10]$, whose exact solution is given by $f(t) = \G^{-1}(1.5) \sqrt{t}, t > 0$. It is shown in Figure \ref{fig_error_amplification} how the relative error of the GL approximation increases artificially above $t \approx 3.4$.

\begin{figure}[ht!]
    \centering
    \includegraphics[width = 0.70\columnwidth]{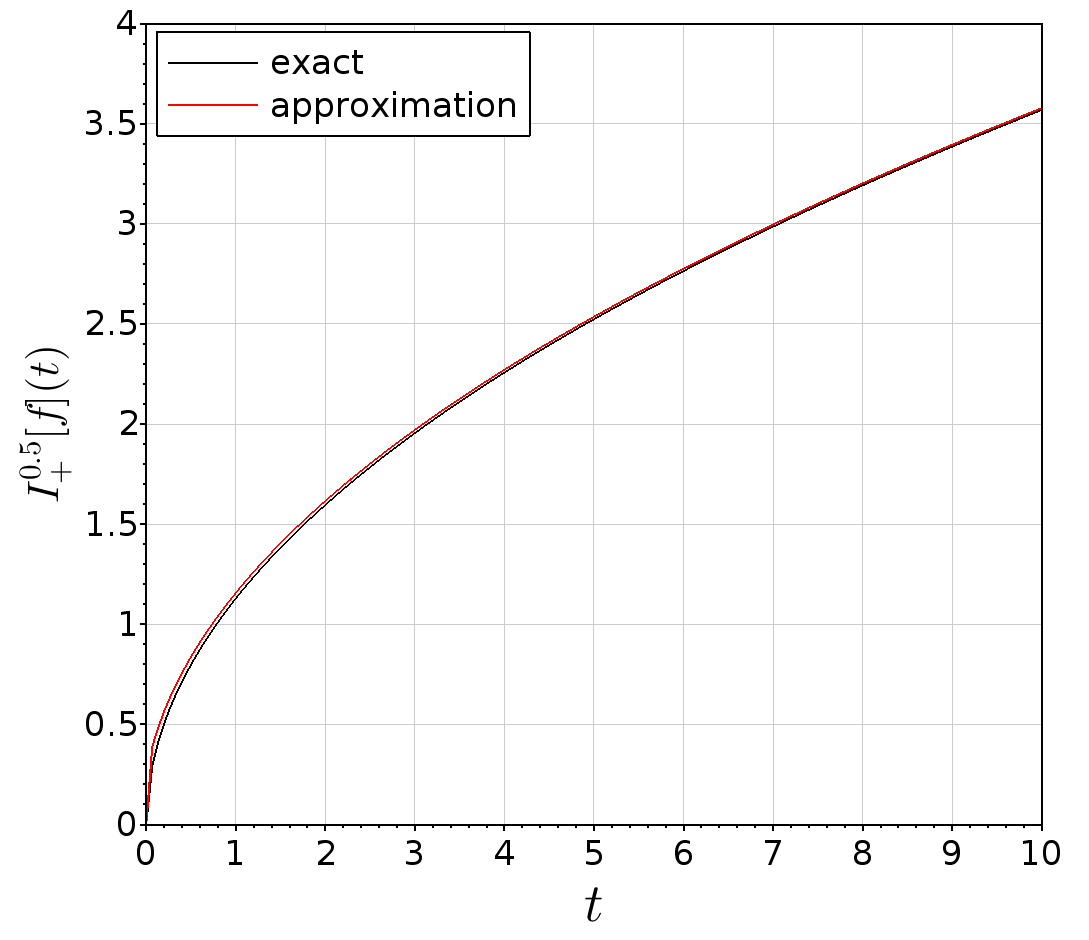}
    \caption{Gr\"unwald-Letnikov approximation with coarser grid $(N = 150)$.}
    \label{fig_coarse_grid}
\end{figure}
\begin{figure}[ht!]
    \centering
    \includegraphics[width = 0.70\columnwidth]{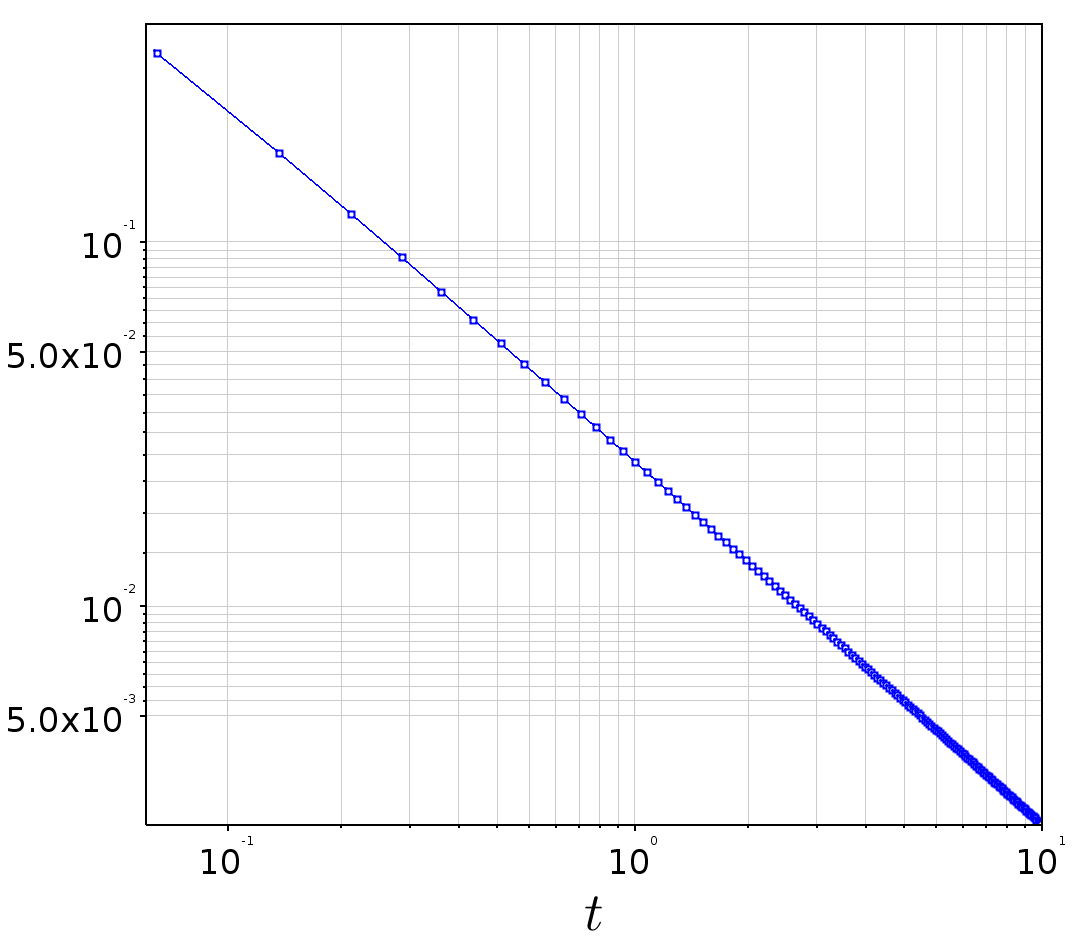}
    \caption{Error in the Gr\"unwald-Letnikov approximation with coarser grid.}
    \label{fig_coarse_grid_error}
\end{figure}

An immediate solution for this issue consists in making the grid coarser; vide Figures \ref{fig_coarse_grid} and \ref{fig_coarse_grid_error}. However, it should be noticed that Euler's scheme, even in the fractional case \cite{lubich_1986}, is an order one method. Therefore, making $\Del t$ large leads to a progressive loss of the accuracy.

The observed errors are due solely to number representation limitations on computers, and have nothing to do with the nature of the fractional operations discussed so far. On the other hand, the rapid decay of GL approximation coefficients introduces some interesting aspects. It shows that, at least for the case of fractional derivative $(\al > 0)$, the lower terminal $a$ (taken equal to zero here) contributions are not always fundamental for the evaluation of $D^{\al}[\cdot]$. In other words, there is some sort of  ``numerically induced locality'' in a  non-local operator. That means that one could use only few past entries for each  $n$ in order to approximate $D^{\al}_n$  and, therefore, $I^{-\al}_n$, according Section  \ref{sec_fundamental_concepts}. This is what Li and Zeng call short-memory principle \cite{li_2015}. In their work, authors digress on how one can estimate a truncation term. We, however, follow a different path, in which zero and one order approximations of the generalized Newton-Cotes are used \cite{li_2015, li_2011}. Such approach provides similar results than those obtained by FLMM methods (up to order one) and does not suffer from the former approach practical complications, such as the need for Laplace transform inversions and truncation terms estimation. 

%% file: numerical_approximation/sec_interpolation_methods.tex
Instead of using FLMM to approximate $I^{\al}_{n}$, one can appeal to interpolation methods, whose formalism is built on orthogonal polynomials local expansions \cite{dahlquist_2003}. The main idea consists in approximating the integrand by an easily integrable function within a given partition of the original interval. As shown by Li and Zheng \cite{li_2015}, this classical method can be extended without further difficulties to the fractional scenario.
For example, let us consider a zero order approximation for $f$, i.e., when one wants to approximate $f(t)$ in a given partition $t \in [t_k, t_{k+1}]$ by $f_{k + 1/2} \equiv 0.5(f_k + f_{k+1})$. In the case of $f(t_k) \approx f_k$ , the left-sided RL integral reads
\begin{eqnarray}\label{eq_NC0_coeff_derivation}
    I_{n}^{\al} &   =       &   \frac{1}{\G(\al)}\int\limits_{0}^{t_n}{(t_n - t')^{\al-1}f(t')dt'}  \nonumber   \\
                &   \approx &   \frac{1}{\G(\al)}\sum\limits_{k = 0}^{n - 1}{\int\limits_{t_k}^{t_{k+1}}{(t_n - t')^{\al-1}f_k dt'}}    \nonumber   \\
                &   \approx &   \frac{1}{\G(\al + 1)}\sum\limits_{k = 0}^{n - 1}{f_k    \left[     (t_n-t_k)^{\al} - (t_n-t_{k+1})^{\al}    \right]}    \nonumber   \\
                &   \approx &   \frac{(\Del t)^{\al}}{\G(\al + 1)}\sum\limits_{k = 0}^{n - 1}{f_k    \left[     (n - k)^{\al} - (n - k - 1)^{\al}    \right]}   \nonumber   \\
                &   \approx &   \sum\limits_{k = 0}^{n - 1}{f_k c_{(n-1)-k}}, 
\end{eqnarray}
in which
\begin{equation}\label{eq_NC0_left_coeff}
    c_k \equiv \frac{(\Del t)^{\al}}{\G(\al + 1)}[(k + 1)^{\al} - k^\al].
\end{equation}
The case for  $f(t_k) \approx f_{k+1}$ is completely analogous, leading to
\begin{equation}\label{eq_NC0_right_approximation}
    I_{n}^{\al} \approx \sum\limits_{k = 0}^{n-1}{f_{k+1}c_{(n-1)-k}},
\end{equation}
where $c_k$ is given by Eq. \eqref{eq_NC0_left_coeff}. Taking the average between Eqs \eqref{eq_NC0_coeff_derivation} and \eqref{eq_NC0_right_approximation}, one deals with the fractional composite trapezoid rule \cite{li_2015}:
\begin{equation}\label{eq_frac_composite_trapezoid}
    I_{n}^{\al} \approx \sum\limits_{k = 0}^{n-1}{f_{k + \frac{1}{2}} c_{(n-1)-k}}.
\end{equation}

In principle, one might choose higher order approximations for $f$ and apply the same procedure above in order to obtain better estimations. This is possible since the integrand can be locally approximated by a polynomial. In the general case, one can represent $f(t)$ by its Lagrange polynomial \cite{ralston_2001, dahlquist_2003}, given by
\begin{equation}\label{eq_lagrange_interpolation}
    f(t) = \sum\limits_{j = 0}^{P-1}{f_{\left(k + \frac{j}{P - 1}\right)} \ell_{\left(k + \frac{j}{P - 1}\right)}(t)}, 
\end{equation}
with
\begin{equation}\label{eq_lagrange_interpolator}
    \ell_{\left(k + \frac{j}{P - 1}\right)}(t) = \prod\limits_{\substack {i = 0 \\ j \neq i}}^{P-1}{\left[   \frac{t- t_{\left(i + \frac{j}{P - 1}\right)}}{t_{\left(k + \frac{j}{P - 1}\right)} - t_{\left(k + \frac{i}{P - 1}\right)}}   \right]}.
\end{equation}
In the expressions above, $P$ is the number of sub-partitions of the interval $[t_k, t_{k+1}]$.

In a formal sense, there is no impediment to the derivation of higher orders schemes, provided by the Fractional Newton-Cotes formula \cite{li_2015}:
\begin{equation}\label{eq_frac_NC}
    \frac{1}{\G(\al)} \sum\limits_{k = 0}^{n-1}{   \sum\limits_{j = 0}^{P-1}{f_{\left(k + \frac{j}{P-1}\right)}   \int\limits_{t_k}^{t_{k+1}}{\ell_{\left(k + \frac{j}{P-1}\right)}}(t')(t_n-t')^{\al-1}dt'      }     }
\end{equation}
It may be tempting to increase  in order to achieve better results; however, one should notice that for equidistant points, Lagrange interpolation presents oscillatory behavior near terminals \cite{dahlquist_2003}, leading to bad representations of the integrand.  In fact, for sufficiently smooth integrands, low order methods provide satisfactory results.

%% file: numerical_approximation/sec_numerical_examples.tex
To illustrate the use of fractional linear multi-step and interpolation methods to approximate left-sided RL integrals, let us consider two integrands:  $f(t) = 1$ and $g(t) = e^t$, whose exact integrals are given respectively by
\begin{equation}\label{eq_exact_solutions}
    \begin{cases}
        I^{\al}[f](t)   =   \frac{1}{\G(\al+1)}t^\al, ~\text{and}   \\
        \\
        I^{\al}[g](t)   =   \frac{1}{\G(\al)}e^t \g(t,\al)     
    \end{cases}
\end{equation}
where $\g(\cdot,\cdot)$ is the lower incomplete gamma function, given by
\begin{equation}\label{eq_lower_incomp_gamma}
    \g(t,\al) = \int\limits_{0}^{t}{u^{\al-1}e^{-u}du},
\end{equation}
whose numerical evaluation can be performed through Taylor series expansion, at least for small $t$ \cite{abergel_2016}:
\begin{equation}\label{eq_series_exp_lower_incomp_gamma}
    \g(t,\al) = \sum\limits_{n = 0}^{\infty}{\frac{(-1)^{n}}{n!(\al+n)}t^{\al + n}}
\end{equation}

It is interesting to notice that both Gr\"unwald-Letnikov approximation and zero order Newton-Cotes method (NC-0) can be represented by a discrete causal convolution:
\begin{equation}\label{eq_discrete_conv}
    I_{n}^{\al} = \sum\limits_{k = 0}^{n}{f_k c_{n-k}}, 
\end{equation}
in which the coefficients are given respectively by
\begin{equation}\label{eq_conv_coeffs}
    \begin{cases}
        c^{\text{GL}}_k = (\Del t)^{\al}(-1)^k \binom{-\al}{k}, ~\text{and}   \\
        \\
        c^{\text{NC-0}}_k = \frac{(\Del t)^{\al}}{\G(\al + 1)}[(k+1)^\al - k^\al]
    \end{cases}
\end{equation}
Therefore, both methods can be efficiently implemented by means of use of FFT algorithms, largely available \cite{scilab_2012}.

\begin{figure}[ht!]
    \centering
    \includegraphics[width = 0.70\columnwidth]{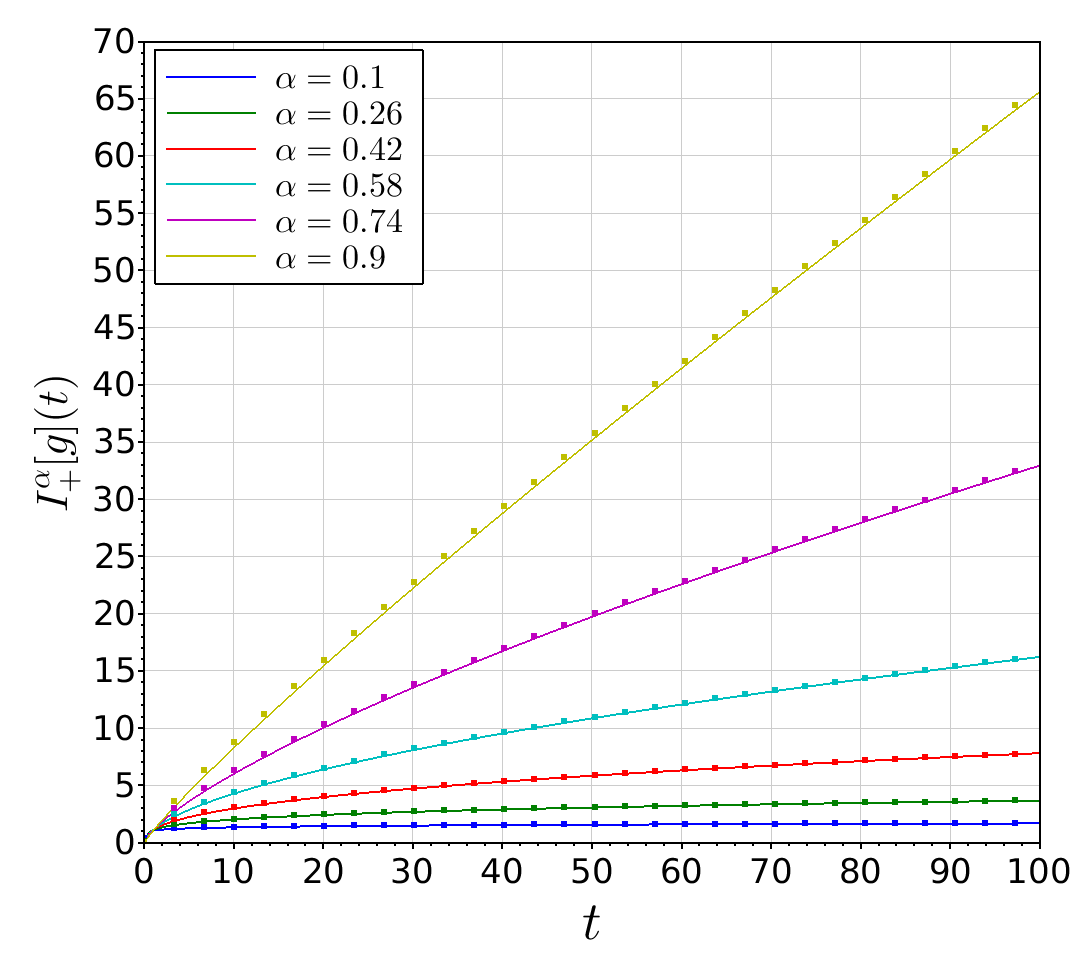}
    \caption{GL approximation for $f(t)$. Lines show the exact result, while approximate solutions are represented by the squares. One should notice how solutions progressively approach the classical solution $f(t) = t$.}
    \label{fig_GL_approximation_for_f}
\end{figure}

\begin{figure}[ht!]
    \centering
    \includegraphics[width = 0.70\columnwidth]{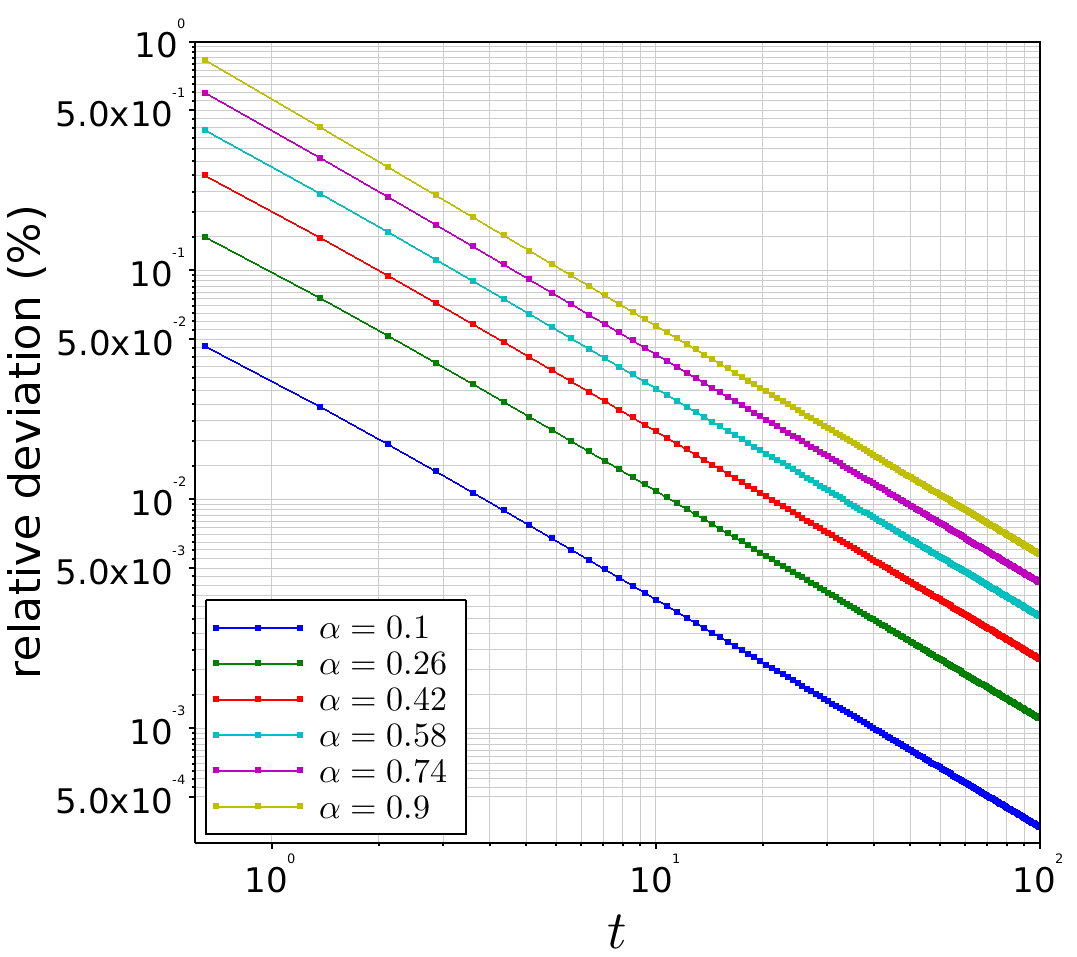}
    \caption{Error of the GL approximation. Due to the singularity at $t = 0$, one observes a larger error near origin.}
    \label{fig_error_GL_approximation_for_f}
\end{figure}

Fig. \ref{fig_GL_approximation_for_f} shows results for $I[f]^{\al}(t)$ using GL approximation, where it can be noticed a good agreement among exact and approximate results. In Fig. \ref{fig_error_GL_approximation_for_f} are shown the errors of the GL approximation, corroborating to the good accuracy statement, since far from origin (a point of singularity), relative deviations fall below 0.1\% for all tested $\al$, in the worst case scenario, at $t \gtrsim 6.0$.
NC-0 method provides similar results to those shown in Fig. \ref{fig_GL_approximation_for_f}. In fact, both methods differ little, as can be seen in Fig. 8, where absolute difference between approximations are shown. In either case, it is important to keep in mind that the integrand  is a zero-order polynomial, which by construction of the methods, must lead to a high fidelity approximation.
\begin{figure}[ht!]
    \centering
    \includegraphics[width = 0.70\columnwidth]{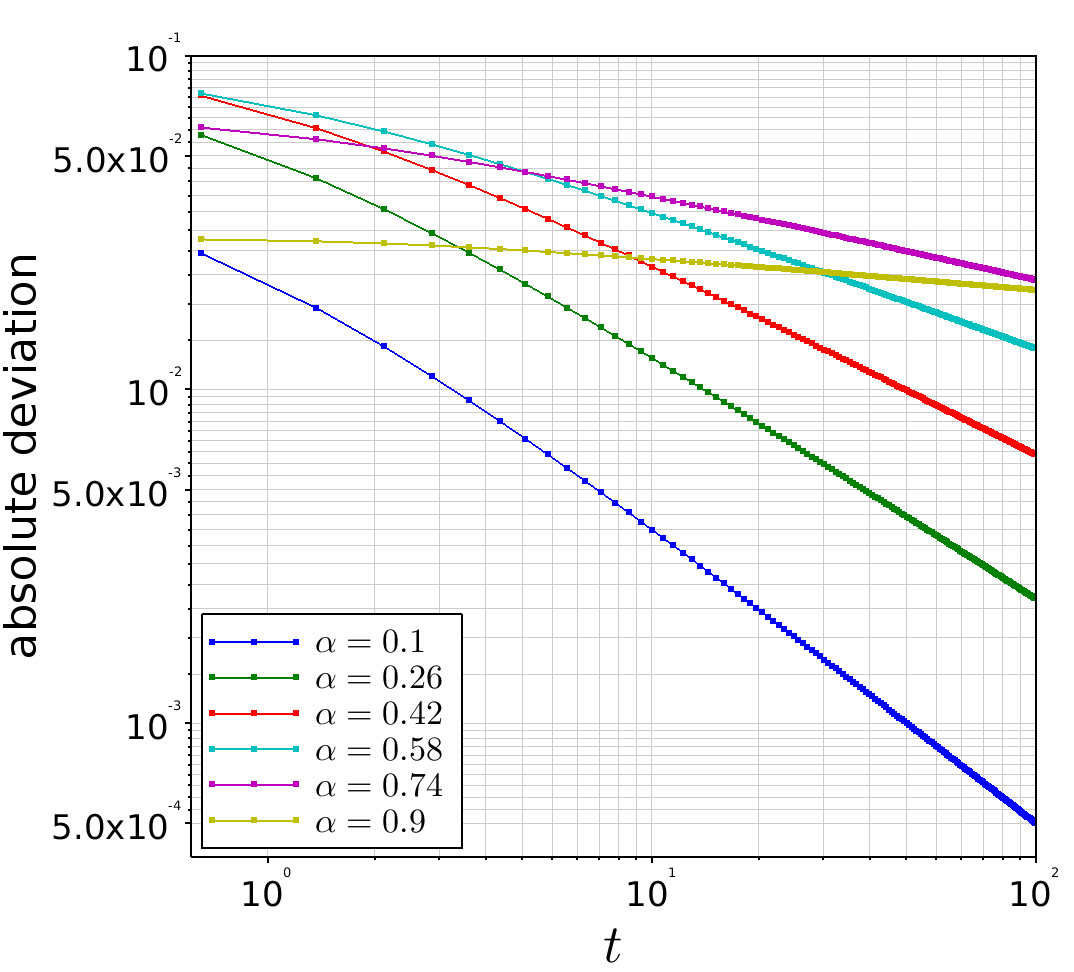}
    \caption{Absolute difference between GL and NC-0 approximations.}
    \label{fig_error_NC0_GL.png}
\end{figure}

\begin{figure}[ht!]
    \centering
    \includegraphics[width = 0.70\columnwidth]{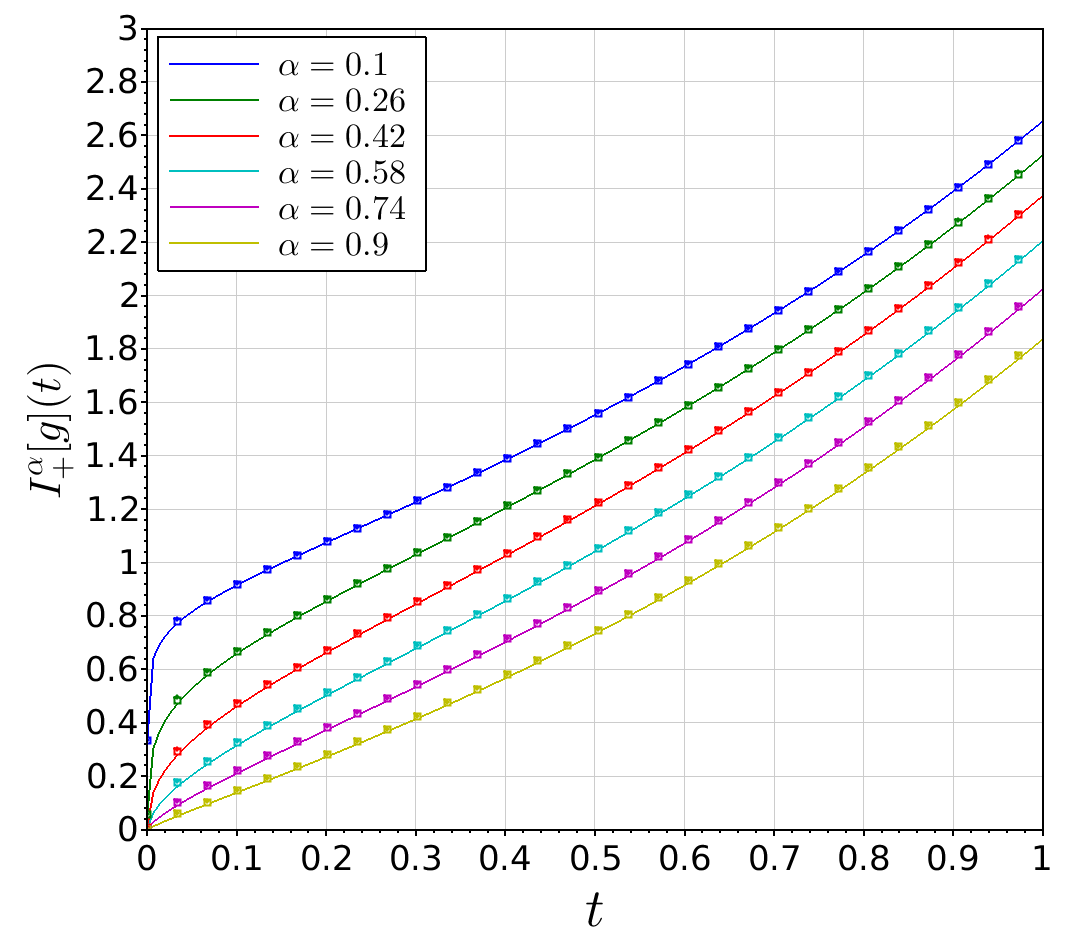}
    \caption{GL approximation for $g(t)$. Lines show the exact result, while approximate solutions are represented by squares.}
    \label{fig_GL_approximation_for_g}
\end{figure}

\begin{figure}[ht!]
    \centering
    \includegraphics[width = 0.70\columnwidth]{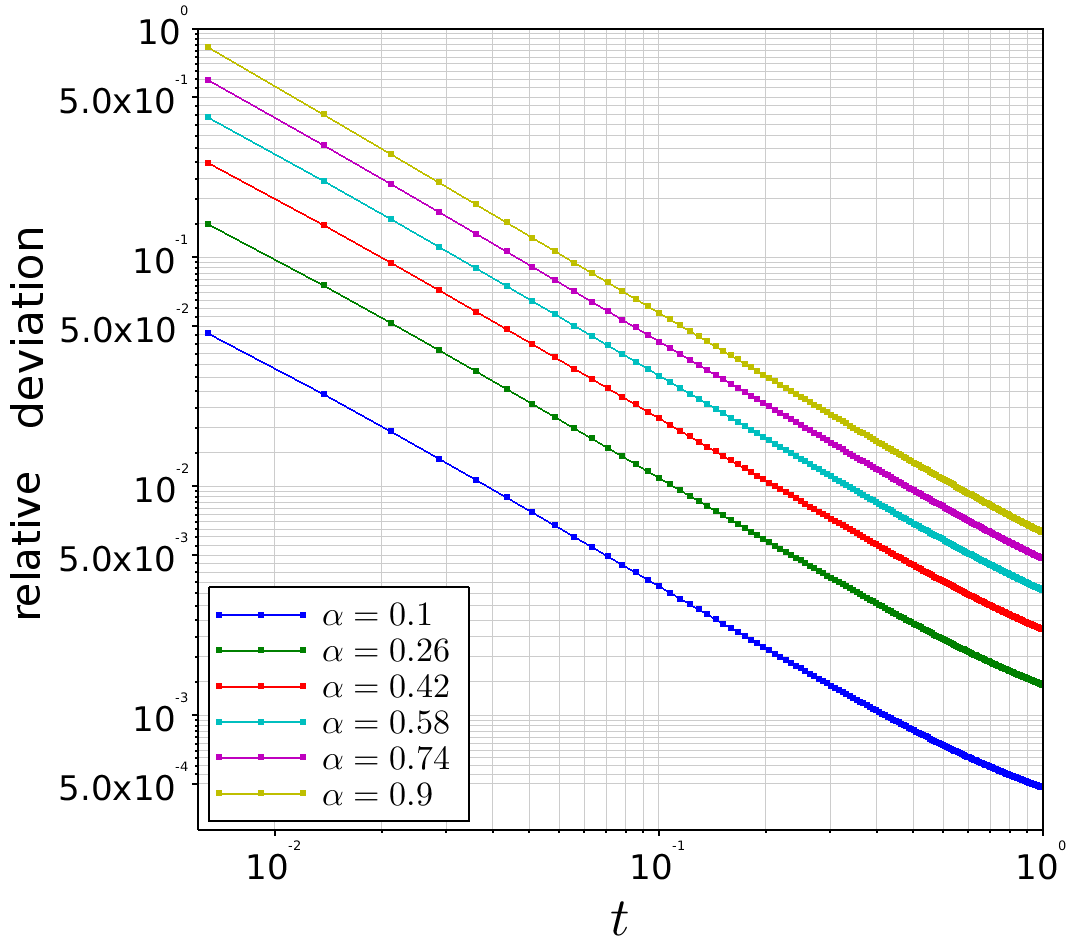}
    \caption{Error in the GL approximation for $g(t)$.}
    \label{fig_error_GL_approximation_for_g}
\end{figure}

\begin{figure}[ht!]
    \centering
    \includegraphics[width = 0.70\columnwidth]{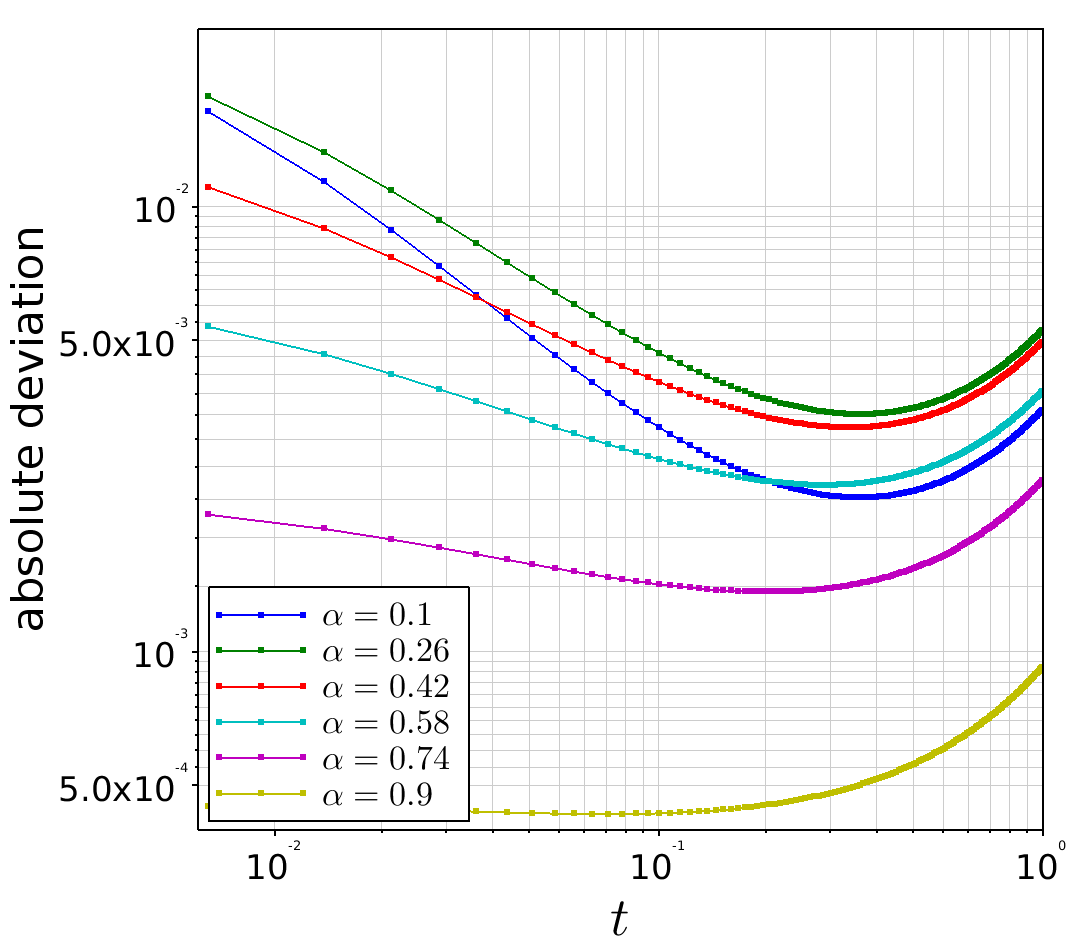}
    \caption{Absolute deviation among GL and NC-0 approximations for $g(t)$. One notices an error amplification for values near superior terminal.}
    \label{fig_error_NC0_GL_for_g}
\end{figure}

Fig. \ref{fig_GL_approximation_for_g} shows results for $I^{\al}[g](t)$ using GL approximation, presenting again a good agreement between exact and approximate solutions. Errors are shown in Fig. \ref{fig_error_GL_approximation_for_g}.
NC-0 method also led to good approximations, but since neither of methods are built to match higher order polynomials, it was verified an increasing of the error magnitude for larger values of $t$, vide Fig. \ref{fig_error_NC0_GL_for_g}. This is also related to the limited extent of the series representation for the lower incomplete gamma function, vide Eqs. \eqref{eq_lower_incomp_gamma} and \eqref{eq_series_exp_lower_incomp_gamma}.

%% file: conclusions/sec_conclusions.tex
The present paper introduced Fractional Calculus main concepts, performing a connection between this theoretical background to applied classical electrodynamics. Such connection is feasible due mainly to the ``convolutional'' structure of both fractional operators -- such Riemman-Liouville integrals and Gr\"unwald-Letnikov derivatives -- and constitutive relations of linear non-local media.

In order to provide means for numerical experimentation on this theoretical development, a review on two of the most simple numerical tools was carried out, in which were pointed out practical aspects, such as gamma function overflow and implementation through FFT algorithms. With respect to the former occurrence, it is clear that order zero Newton-Cotes approximation should be chosen when fine grids are necessary.

Despite the fact that only low order schemes were used in this work, general formulations were exposed, allowing one to build more complex models, at the expense of either performing a Laplace transform inversion (FLMM), or using higher order polynomial approximations of the integrand (fractional Newton-Cotes method). Both choices have their virtues and drawbacks and it is on the operator to decide which one is better for a given problem.

%% file: acknowledgements/sec_acknowledgements.tex
The author would like to thank D.Sc. Luis A. R. Ramirez from Military Institute of Engineering -- IME at Rio de Janeiro and Helen Braga, from CBPF for their support and critical analysis of this work.